# LAGRANGIAN MULTIPLIERS
# FOR PRESUBCONVEXLIKE OPTIMIZATION PROBLEMS
# OF SET-VALUED MAPS

Renying Zeng[1]
School of Mathematical Sciences, Chongqing Normal University, China

## ABSTRACT

In this paper, we discuss scalar Lagrangian multipliers and vector Lagrangian multipliers for constrained set-valued optimization problems. We obtain some necessary conditions, sufficient conditions, as well as necessary and sufficient conditions for the existence of weakly efficient solutions.

[1] Mailing address: Mathematics Department, Saskatchewan Polytechnic, 1130 Idylwyld Dr. N, Saskatoon, Saskatchewan, Canada S7L 4J7. Email: renying.zeng@saskpolytech.ca

# 1. Introduction

Many investigations have been made for convex or generalized convex multi-objective decision or optimization problems. For example, P. L. Yu [17] discussed the non-dominated solutions in cone-convex multi-objective decision problems; J. Borwein [8] established a general alternative theorem for convex set-valued maps in a real topological linear space and got a multiplier principle; O. Ferrero [11] obtained a alternative theorem and related results for a set-valued optimization problem; H. W. Corley [13] established a Lagrangian duality theorem and some optimality conditions for vector optimization with set-valued maps. We notice that, recently, more authors have paid their attentions to vector optimization problems with cone convex and set-valued maps, as well as their applications in other disciplines, e.g., see Refs. [3, 4, 9-13, 22, 25, 26, 28, 30, 31, 35, 37, 39, 40].

In this paper, we establish some necessary conditions, sufficient conditions, as well as necessary and sufficient conditions for the existence of a weakly efficient solution for a set-valued minimization problem. Our necessary conditions, sufficient conditions, as well as necessary and sufficient conditions are related to Lagrangian multipliers of the minimization problem, which has not only equality constraints but also inequality constraints.

# 2. Preliminaries

Throughout this paper, scalars of linear space are always real. A topological linear space $Y$ with a convex cone $Y_+$ is said to be an ordered topological linear space, the partial order in $Y$ is induced by:

$$y^1 \geq y^2, \text{ iff } y^1 - y^2 \in Y_+,$$

$$y^1 > y^2, \text{ iff } y^1 - y^2 \in \text{int } Y_+,$$

where $\text{int } Y_+$ is the interior of $Y_+$, which is nonempty, i.e., $\text{int } Y_+ \neq \emptyset$.

Suppose that $Y^*$ is the topological dual of $Y$, then

$$Y_+^* = \{\xi \in Y^* : \xi(y) \geq 0, \forall y \in Y_+\}$$

is said to be the dual cone of $Y_+$.

A convex cone $Y_+$ of $Y$, with apex at the origin $0_Y$ of $Y$, is said to be a positive cone of $Y$.

A convex cone $Y_+$ of $Y$ is said to be pointed if $Y_+ \cap (-Y_+) = \{0_Y\}$.

In this paper, we always assume that every convex cone is positive, closed, pointed and with nonempty interior.

Suppose $f: X \to 2^Y$ is a set-valued function, where $X$ is a linear space, $2^Y$ denotes the power set of $Y$. For any nonempty $D \subseteq X$, Setting

$$f(D) = \bigcup_{x \in D} f(x).$$

For $x \in X, \xi \in Y^*$, we may set

$$\xi(f(x)) \geq 0, \text{ iff } \xi(y) \geq 0, \forall y \in f(x);$$

$$\xi(f(D)) \geq 0, \text{ iff } \xi(f(x)) \geq 0, \forall x \in D.$$

We denote by $R$ the set of real numbers. For $A, B \subseteq R$, write

$$A \geq B, \text{ if } a \geq b, \forall a \in A, \forall b \in B. \tag{2.1}$$

**Definition 2.1** A set-valued map $u: D \to 2^Y$ is said to be bounded if $\bigcup_{x \in D} u(x)$ is a bounded subset of $Y$.

**Lemma 2.1** [6] Let $C, D \subseteq Y$ be two nonempty convex sets. If $\text{int } C \neq \varnothing, D \cap \text{int } C \neq \varnothing$ then

$$\exists 0_Y \neq \xi \in Y^*: \xi(\text{int } C) \subseteq \xi(D) + \text{int } Y_+, \xi(C) \subseteq \xi(D) + Y_+.$$

**Lemma 2.2** [7] If $\xi \in Y_+^* \setminus \{0_{Y^*}\}$, then $\xi(y) > 0, \forall y \in \text{int } Y_+$.

Suppose that $X$ is a nonempty set and $D$ is a nonempty subset of $X$, and $Y$ is an ordered topological linear space with a convex cone $Y_+$.

**Definition 2.2** A set-valued map $f: D \to 2^Y$ is said to be convexlike on $D$ if $\forall x^1, x^2 \in D, \forall \alpha \in (0,1)$

$$\alpha f(x^1) + (1-\alpha) f(x^2) \subseteq f(D) + Y_+.$$

**Definition 2.3** A set-valued map $f: D \to 2^Y$ is said to be subconvexlike on $D$ if $\exists$ bounded function $u: X \to 2^Y$ such that $\forall x^1, x^2 \in D, \forall \alpha \in (0,1), \forall \varepsilon > 0$

$$\varepsilon u + \alpha f(x^1) + (1-\alpha) f(x^2) \subseteq f(D) + Y_+.$$

**Definition 2.4** A set-valued map $f: D \to 2^Y$ is said to be presubconvexlike on $D$ if $\exists$ bounded function $u: X \to 2^Y$ such that $\forall x^1, x^2 \in D, \forall \alpha \in (0,1), \forall \varepsilon > 0, \exists \tau > 0$

$$\varepsilon u + \alpha f(x^1) + (1-\alpha) f(x^2) \subseteq \tau f(D) + Y_+ .$$

**Theorem 2.1** (i) The set-valued map $f: X \to 2^Y$ is convexlike on $D$ if and only if $f(D) + Y_+$ is convex;

(ii) The map set-valued map $f: X \to 2^Y$ is subconvexlike on $D$ if and only if $f(D) + \text{int } Y_+$ is convex;

(iii) The map set-valued map $f: X \to 2^Y$ is presubconvexlike on $D$ if and only if $\bigcup_{t>0}(tf(D) + \text{int } Y_+)$ is convex, where $t > 0$ is any positive scalar.

**Proof.** We refer to Paeck [5, Lemma 2.2] for the proof of ( i) and ( ii). We are only going to show (iii).

Let $c_1, c_2 \in \bigcup_{t>0}(tf(D) + \text{int } Y_+)$, $\alpha \in (0,1)$. Then there exist $x^1, x^2 \in D, y_1, y_2 \in \text{int } Y_+$, and positive numbers $t_1, t_2$ such that

$$c^1 \in t_1 f(x^1) + y_1, c^2 \in t_2 f(x^2) + y_2.$$

Since $f$ is presubconvexlike on $D$, there exists a bounded set-valued map $u: D \to 2^Y$, $x^3 \in D$, and a positive number $\tau$ such that $\forall y^1 \in f(x^1), \forall y^2 \in f(x^2)$

$$\varepsilon u + \frac{\alpha t_1}{\alpha t_1 + (1-\alpha)t_2} y^1 + \frac{(1-\alpha)t_2}{\alpha t_1 + (1-\alpha)t_2} y^2 \subseteq tf(x^3) + Y_+ \tag{2.2}$$

Because $\text{int } Y_+$ is a convex cone, we have $y_0 := \alpha y_1 + (1-\alpha)y_2 \in \text{int } Y_+$. It deduces that there exists a neighborhood $U$ of the origin of $Y$ such that $y_0 + U \subseteq \text{int } Y_+$. Since $U$ is absorbent and $u$ is bounded, we may choose a positive number $\varepsilon$ which is small enough such that $-(\alpha t_1 + (1-\alpha)t_2)\varepsilon u \in U$. And so

$$y_0 - (\alpha t_1 + (1-\alpha)t_2)\varepsilon u \in y_0 + U \subseteq \text{int } Y_+. \tag{2.3}$$

Therefore combining (2.2) and (2.3)

$$\alpha c^1 + (1-\alpha)c^2$$

$$= \alpha t_1 y^1 + (1-\alpha)t_2 y^2 + \alpha y_1 + (1-\alpha)y_2$$

$$= (\alpha t_1 + (1-\alpha)t_2)[\frac{\alpha t_1}{\alpha t_1 + (1-\alpha)t_2} y^1 + \frac{(1-\alpha)t_2}{\alpha t_1 + (1-\alpha)t_2} y^2] + y_0$$

$$\subseteq (\alpha t_1 + (1-\alpha)t_2)(tf(x^3) + Y_+ - \varepsilon u) + y_0$$

$$= (\alpha t_1 + (1-\alpha)t_2) tf(x^3) + (\alpha t_1 + (1-\alpha)t_2) Y_+ - (\alpha t_1 + (1-\alpha)t_2)\varepsilon u + y_0$$

$$\subseteq (\alpha t_1 + (1-\alpha)t_2) tf(D) + Y_+ + \text{int } Y_+$$

$$\subseteq (\alpha t_1 + (1-\alpha)t_2) \tau f(D) + \text{int } Y_+ \subseteq \bigcup_{t>0}(tf(D) + \text{int } Y_+)$$

Therefore $\bigcup_{t>0}(tf(D) + \text{int } Y_+)$ is convex.

On the other hand, assume $\bigcup_{t>0}(tf(D) + \text{int } Y_+)$ is a convex set.

Let $y^0 \in \text{int } Y_+, x^1, x^2 \in D, \alpha \in (0,1), \varepsilon > 0$, and $y^1 \in f(x^1), y^2 \in f(x^2)$. Since

$$f(x^i) + y^0 \subseteq \bigcup_{t>0}(tf(D) + \text{int } Y^+), i = 1,2; \text{ and } \bigcup_{t>0}(tf(D) + \text{int } Y_+) \text{ is convex, we get}$$

$$\varepsilon y^0 + \alpha y^1 + (1-\alpha)y^2 = \alpha(y^1 + \varepsilon y^0) + (1-\alpha)(y^2 + \varepsilon y^0)$$

$$\subseteq \bigcup_{t>0}(tf(D) + \text{int } Y^+)$$

Therefore there exists positive number $t > 0$, such that

$$\varepsilon u + \alpha f(x^1) + (1-\alpha)f(x^2) \subseteq tf(D) + \text{int } Y_+ \subseteq tf(D) + Y_+.$$

Where $u \equiv \{y^0\} : D \to 2^Y$ is a bounded set-valued map.

Hence $f$ is presubconvexlike. The proof is complete. □

From the definitions we see that: convexlike $\Rightarrow$ subconvexlike $\Rightarrow$ presubconvexlike.

Since any real-valued function is a convexlike function, it is clear that subconvexlikeness does not imply convexlikeness. However, the following example illustrates that presubconvexlikeness does not imply subconvexlikeness, either.

**Example 2.1** Let

$$D = \{(x_1, x_2) \in R^2 : x_1 \geq 0, x_2 \geq 0, x_1^2 + x_2^2 \geq 1\},$$

and

$$Y = R^2, Y_+ = R_+^2 = \{(x_1, x_2) \in R^2 : x_1 \geq 0, x_2 \geq 0\},$$

where $R^2$ is the 2-dimensional Euclidean space.

Define a set-valued map $f: D \to 2^Y$ by

$$f(x_1, x_2) = \{(x_1, x_2)\} \cup \{(a,b): a \geq 0, b \geq 0, a^2 + b^2 = 1\}.$$

We note that

$$f(D) + \text{int } Y_+$$

$$= \{(x_1, x_2): x_1 > 0, x_2 > 0, x_1^2 + x_2^2 > 1\}$$

is not a convex set, so $f$ is not subconvexlike. But

$$\bigcup_{t>0}(tf(D) + \text{int } Y_+) = \text{int } Y_+$$

is convex, therefore $f$ is presubconvexlike. □

**Theorem 2.2** Assume that $f, g, h$ satisfy

(a1) There exist bounded set-valued maps $u_1: X \to 2^Y, u_2: X \to 2^Z$, for which

$\forall \alpha \in (0,1), \forall x^i \in D, \forall y^i \in f(x^i), \forall z^i \in g(x^i), \forall w^i \in h(x^i), i = 1,2, \forall \varepsilon > 0, \exists \tau_i > 0, i = 1, 2, 3,$

such that

(a1.1) $\quad \varepsilon u_1 + \alpha y^1 + (1-\alpha)y^2 \in \tau_1 f(D) + Y_+,$

(a1.2) $\quad \varepsilon u_2 + \alpha z^1 + (1-\alpha)z^2 \in \tau_2 g(D) + Z_+,$

(a1.3) $\quad \alpha w^1 + (1-\alpha)w^2 \in \tau_3 h(D);$

(a2) $\text{int } h(D) \neq \emptyset;$

and ( i) and ( ii) denote the systems

(i) $\exists x \in D, s.t., f(x) \cap (-\text{int } Y_+) \neq \emptyset, g(x) \cap (-Z_+) \neq \emptyset, 0_W \in h(x);$

(ii) $\exists (\xi, \eta, \zeta) \in (Y_+^* \times Z_+^* \times W^*) \setminus \{(0_Y, 0_Z, 0_W)\}$ such that

$$\xi(f(x)) + \eta(g(x)) + \varsigma(h(x)) \geq 0, \ \forall x \in D.$$

If (i) has no solution then (ii) has solutions.

Moreover if (ii) has a solution $(\xi, \eta, \varsigma)$ with $\xi \neq 0_{Y^*}$ then (i) has no solutions.

**Lemma 2.3** The assumption (a1) is satisfied if and only if the following set $B = \subseteq Y \times Z \times W$ is convex:

$$B = \{(y,z,w): y \in \bigcup_{t>0}(tf(D) + \text{int } Y_+), z \in \bigcup_{t>0}(tg(D) + \text{int } Z_+), w \in \bigcup_{t>0}th(D)\}$$

$$= (\bigcup_{t>0}(tf(D) + \text{int } Y_+)) \times (\bigcup_{t>0}(tg(D) + \text{int } Z_+)) \times (\bigcup_{t>0}th(D)).$$

**Proof.** ($\Rightarrow$) Suppose that (a1) is satisfied.

Similar to Lemma 2.3 $\bigcup_{t>0}(tf(D) + \text{int } Y_+)$, $\bigcup_{t>0}(tg(D) + \text{int } Z_+)$ is convex.

Let $\overline{w}^i \in \bigcup_{t>0}th(D), i = 1,2$, then $\exists x^i \in D$, $\exists$ positive numbers $t_i > 0$ such that $\overline{w}^i = t_i h(x^i), i = 1,2.$ $\exists w^i \in h(x^i)$,

$$\overline{w}^i = t_i w^i, i = 1,2.$$

By the assumption (a1), $\forall \alpha \in (0,1)$, $\exists x^3 \in D, \tau_3 > 0$ such that

$$\alpha t_1 \overline{w}^1 + (1-\alpha)t_2 \overline{w}^2$$

$$= (\alpha t_1 + (1-\alpha)t_2)[\frac{\alpha t_1}{\alpha t_1 + (1-\alpha)t_2}w^1 + \frac{(1-\alpha)t_2}{\alpha t_1 + (1-\alpha)t_2}w^2]$$

$$\subseteq (\alpha t_1 + (1-\alpha)t_2)\tau h(x^3)$$

$$\subseteq \bigcup_{t>0} th(D)$$

Hence $\bigcup_{t>0} th(D)$ is convex.

So, $B = (\bigcup_{t>0}(tf(D) + \text{int } Y_+)) \times (\bigcup_{t>0}(tg(D) + \text{int } Z_+)) \times (\bigcup_{t>0} th(D))$ is convex.

($\Leftarrow$) Assume $B$ is a convex set. Then $\bigcup_{t>0}(tf(D) + \text{int } Y_+), \bigcup_{t>0}(tg(D) + \text{int } Z_+)$, and $\bigcup_{t>0} th(D)$ are all convex sets. Again, from Lemma 2.3, (a1.1) and (a1.2) are satisfied.

Since $\bigcup_{t>0} th(D)$ is convex, $\forall \alpha \in (0,1), \forall x^i \in D, \forall w^i \in h(x^i), i = 1,2,$ we have

$$\alpha w^1 + (1-\alpha)w^2 \in \bigcup_{t>0} th(D).$$

Therefore $\exists \tau_3 > 0$ such that

$$\alpha w^1 + (1-\alpha)w^2 \in \tau_3 h(x).$$

And so (a1.3) is satisfied. The proof of Lemma 2.3 is complete. □

**Proof of Theorem 2.2**.

By Lemma 2.3, $B$ defined above is convex. From (a2), $\text{int } B \neq \emptyset$. We also have $(0_Y, 0_Z, 0_W) \notin B$ since (i) has no solution.

Therefore, according to the separation theorem of convex sets of topological linear space (Lemma 2.1), there exists a nonzero vector $(\xi, \eta, \varsigma) \in Y^* \times Z^* \times W^*$ such that

$$\xi(t_1 y + y^0) + \eta(t_2 z + z^0) + \varsigma(t_3 w) \geq 0 \qquad (2.4)$$

for all $x \in D, y \in f(x), y^0 \in \text{int } Y_+, z \in g(x), z^0 \in \text{int } Z_+, w \in h(x), t_i > 0, i = 1,2,3$.

Since $\text{int } Y_+, \text{int } Z_+$ are convex cones, we get

$$\xi(t_1 y + s_1 y^0) + \eta(t_2 z + s_2 z^0) + \varsigma(t_3 w) \geq 0 \qquad (2.5)$$

for all $x \in D, y \in f(x), y^0 \in \text{int } Y_+, z \in g(x), z^0 \in \text{int } Z_+, w \in h(x), t_i > 0, s_i > 0, i = 1,2$.

Let $t_i \to 0, i = 1,2,3, s_2 \to 0$ in (2.5) we see that

$$\xi(y^0) \geq 0, \forall y^0 \in \text{int } Y_+.$$

Therefore $\xi(y) \geq 0, \forall y \in Y_+$. Hence $\xi \in Y_+^*$.

Similarly, let $t_i \to 0, i = 1,2,3, s_1 \to 0$ in (2.5) we get $\eta \in Z_+^*$. Thus

$$(\xi, \eta, \varsigma) \in Y_+^* \times Z_+^* \times W^*.$$

Let $s_i \to 0, i = 1,2$ in (2.5) again we attain

$$\xi(y) + \eta(z) + \varsigma(w) \geq 0$$

for $\forall x \in D, \forall y \in f(x), \forall z \in g(x), \forall w \in h(x)$. Which implies that (ii) has solutions.

On the other hand, suppose that (ii) has a solution $(\xi, \eta, \varsigma)$ with $\xi \neq 0_{Y_+^*}$. If the system (ii) had a solution $x \in D$, there would exist $y \in f(x), z \in g(x), w \in h(x)$ such that

$$y \in \text{int } Y_+, z \in -Z_+, w = 0_W.$$

Therefore by Lemma 2.2

$$\xi(y) + \eta(z) + \varsigma(w) < 0.$$

Contradicting our assumption (i). Therefore, the system (ii) had no solution. We complete the proof. □

We note that (a2) states $\text{int } B \neq \emptyset$. Actually, to use the separation theorem of two convex sets of a finite dimensional space, we do not need the condition $\text{int } B \neq \emptyset$. Therefore, by the proof of Theorem 2.2, we get following Corollary 2.1.

**Corollary 2.1** Let $Y$, $Z$, and $W$ are finite dimensional. Assume that $\exists$ bounded set-valued maps $u_1: X \to 2^Y, u_2: X \to 2^Z$, and $\forall \alpha \in (0,1), \forall x^i \in D, \forall y^i \in f(x^i)$,

$\forall z^i \in g(x^i), \forall w^i \in h(x^i), i = 1,2, \forall \varepsilon > 0, \exists \tau_i > 0,\ i = 1,2,3$, such that

$$\varepsilon u_1 + \alpha y^1 + (1-\alpha) y^2 \in \tau_1 f(D) + Y_+,$$

$$\varepsilon u_2 + \alpha z^1 + (1-\alpha) z^2 \in \tau_2 g(D) + Z_+,$$

$$\alpha w^1 + (1-\alpha) w^2 \in \tau_3 h(D);$$

and (i) and (ii) denote the systems

(i) $\exists x \in D, s.t., f(x) \cap (-\text{int } Y_+) \neq \emptyset,\ g(x) \cap (-Z_+) \neq \emptyset, 0_W \in h(x)$;

(ii) $\exists (\xi, \eta, \varsigma) \in (Y_+^* \times Z_+^* \times W^*) \setminus \{(0_Y, 0_Z, 0_W)\}$ such that

$$\xi(f(x)) + \eta(g(x)) + \varsigma(h(x)) \geq 0,\ \forall x \in D.$$

If (i) has no solution then (ii) has solutions.

Moreover if (ii) has a solution $(\xi, \eta, \varsigma)$ with $\xi \neq 0_{Y^*}$ then (i) has no solution. □

**Remark 2.1**. For Theorem 2.1 and Corollary 2.1, we notice that if the Nontrivial Abnormal Multiplier Constraint Qualification (NNAMCQ) (Definition 3.3) or the Slater Constraint Qualification (SCQ) (Definition 3.2) is satisfied, then the System (ii) can be read as: $\exists (\xi, \eta, \varsigma)$ with $\xi \neq 0_{Y^*}$ such that

$$\xi(f(x)) + \eta(g(x)) + \varsigma(h(x)) \geq 0, \forall x \in D.$$

## 3. Lagrangian Multipliers

Consider the following optimization problem with set-valued maps:

$$\text{(VP)} \quad \begin{array}{l} Y_+ - \min f(x), \\ \text{s.t., } g(x) \cap (-Z_+) \neq \emptyset, 0_Z \in h(x), \\ x \in X. \end{array}$$

Let $D$ be the feasible set of (VP), i.e.,

$$D = \{x \in X : g(x) \cap (-Y_+) \neq \emptyset, 0_W \in h(x)\}.$$

**Definition 3.1** $\bar{x} \in D$ is said to be a weakly efficient solution of (VP), if $\exists \bar{y} \in f(\bar{x})$ such that $\forall x \in D$, there is no $y \in f(x)$ for which $(\bar{y} - y) \in \text{int } Y_+$.

**Lemma 3.1** $\bar{x} \in D$ is a weakly efficient solution of (VP), if and only if $\exists \bar{y} \in f(\bar{x})$ such that $(\bar{y} - f(D)) \cap \text{int } Y_+ = \emptyset$.

**Definition 3.2** The problem (VP) satisfies the Slater Constraint Qualification (SCQ) if $\forall (\eta, \varsigma) \in (Z_+^* \times W^*) \setminus \{(0_{Z^*}, 0_{W^*})\}$, $\exists x \in D$, and there exists a negative real number $t \in \eta(g(x)) \cap \varsigma(h(x))$.

Our Slater Constraint Qualification is similar to "the regular hypothesis" in V. Jeyakumar [1], and M. S. Bazaraa [16].

**Definition 3.3** Let $\bar{x} \in D$. We say No Nontrivial Abnormal Multiplier Constraint Qualification (NNAMCQ) holds at $\bar{x}$ if

$$\min_{x \in D}[\eta(g(x)) + \varsigma(h(x))] = 0$$
$$\min_{x \in D} \eta(g(\bar{x})) = 0$$
$$\eta \in Z_+^*, \eta \in W^*$$

implies that $(\eta, \varsigma) = 0_{Z^* \times W^*}$.

We note that (NNAMCQ) is weaker than (SCQ).

**Theorem 3.1** Let $\bar{x} \in D$. Assume that $f(x) - f(\bar{x}), g(x), h(x)$ satisfy the generalized convexity condition (a1) and the inner point condition (a2). Then $\bar{x}$ is a weakly efficient solution of (VP) implies $\exists \bar{y} \in f(\bar{x})$ and $\exists$Lagrangian multiplier $\xi \in Y_+^*, \eta \in Z_+^*, \varsigma \in W^*$ with $\xi \neq 0_{Y^*}$ such that

$$\min_{x \in D}[\xi(f(x)) + \eta(g(x)) + \varsigma(h(x))] \geq \xi(\bar{y})$$
$$\min \eta(g(\bar{x})) = 0,$$

where, for the second equality, the minimum is taken for all $z \in g(\bar{x})$.

Inversely, if for $\bar{x} \in D$, (NNAMCQ) holds at $\bar{x}$, and if $\exists$Lagrangian multiplier $\xi \in Y_+^*, \eta \in Z_+^*, \varsigma \in W^*$ with $\xi \neq 0_{Y^*}$ such that

$$\min_{x \in D}[\xi(f(x)) + \eta(g(x)) + \varsigma(h(x))] \geq \xi(\bar{y})$$
$$\min \eta(g(\bar{x})) = 0,$$

then $\bar{x}$ is a weakly efficient solution of (VP).

**Proof**. If $\bar{x}$ is a weakly efficient solution of (VP), then, according to Theorem 2.2, $\exists \bar{y} \in f(\bar{x})$ such that the following system

$$(f(x) - \bar{y}) \cap (-\operatorname{int} Y_+) \neq \varnothing, g(x) \cap (-Z_+) \neq \varnothing, 0_W \in h(x)$$

has no solution for $x \in D$. Hence by Theorem 2.2 $\exists \xi \in Y_+^*, \eta \in Z_+^*, \varsigma \in W^*$ with $(\xi, \eta, \varsigma) \neq 0_{Y^* \times Z^* \times W^*}$ such that

$$\xi(f(x) - \bar{y}) + \eta(g(x)) + \varsigma(h(x)) \geq 0, \forall x \in D.$$

i.e.,

$$\xi(f(x)) + \eta(g(x)) + \varsigma(h(x)) \geq \xi(\bar{y}), \forall x \in D. \tag{3.1}$$

Since $\bar{x} \in D$, $\exists \bar{z} \in g(\bar{x})$ such that $\bar{z} \in -Z_+$. This and $\eta \in Z_+^*$ yield

$$\eta(\bar{z}) \leq 0.$$

On the other hand, due to $0_W \in h(\bar{x})$, taking $x = \bar{x}$ in (3.1) we have

$$\xi(\bar{y}) + \eta(\bar{z}) \geq \xi(\bar{y}).$$

Which followed by

$$\eta(\bar{z}) \geq 0.$$

Therefore

$$\eta(\bar{z}) = 0. \tag{3.2}$$

Taking $x = \bar{x}$ in (3.1) again we attain

$$\xi(\bar{y}) + \eta(g(\bar{x})) + \varsigma(h(\bar{x}))] \geq \xi(\bar{y}). \tag{3.3}$$

Hence

$$\eta(g(\bar{x})) \geq 0.$$

This and (3.2) yield

$$\min \eta(g(\bar{x})) = 0.$$

So, we obtain

$$\xi(\bar{y}) = \xi(\bar{y}) + \eta(\bar{z}) + \varsigma(0_W)$$
$$\in \xi(f(\bar{x})) + \eta(g(\bar{x})) + \varsigma(h(\bar{x})).$$

Therefore, combine (3.2) and (3.3) one has

$$\min_{x \in D}[\xi(f(x)) + \eta(g(x)) + \varsigma(h(x))] = \xi(\bar{y}).$$

Inversely, if for $\bar{x} \in D$, $\exists$Lagrangian multiplier $\xi \in Y_+^*, \eta \in Z_+^*, \varsigma \in W^*$ with $\xi \neq 0_{Y^*}$ such that

$$\min_{x \in D}[\xi(f(x)) + \eta(g(x)) + \varsigma(h(x))] \geq \xi(\bar{y})$$
$$\min \eta(g(\bar{x})) = 0,$$

then

$$\xi(f(x) - \bar{y}) + \eta(g(x)) + \varsigma(h(x)) \geq 0, \forall x \in D. \tag{3.4}$$

If $\xi = 0_{Y^*}$, then (3.4) implies that

$$\eta(g(x)) + \varsigma(h(x)) \geq 0, \forall x \in D. \tag{3.5}$$

Therefore, due to $0_W \in h(x)$ we have

$$\eta(g(x)) \geq 0, \forall x \in D.$$

Hence, $x \in D$ implies that $\exists z \in g(x))$ such that

$$\eta(z) = 0. \tag{3.6}$$

Therefore,

$$\min \eta(g(x)) + \varsigma(h(x)) = 0, \forall x \in D.$$

By the (NNAMCQ) condition, we have

$$(\eta, \varsigma) = 0_{Z^* \times W^*}$$

Which contradicts to the Assumption. Therefore we must have

$$\xi \neq 0_{Y^*}.$$

Now, assume that $\bar{x} \in D$ is not a weakly efficient solution of (VP). Then, by Definition 3.1, $\exists x_0 \in D$ such that $\exists y_0 \in f(x_0)$ for which $\bar{y} - y_0 \in \text{int } Y_+$. So, from Lemma 2.2 we get

$$\xi(y_0 - \bar{y}) < 0. \tag{3.7}$$

Since $x_0 \in D$, $\exists z_0 \in g(x_0)$ such that $z_0 \in -Z_+$. Hence

$$\eta(z_0) \leq 0.$$

This and (3.7) together deduce that

$$\xi(y_0 - \bar{y}) + \eta(g(z_0)) + \varsigma(0_W) < 0.$$

Which is contradicting to (3.4) because $x_0 \in D$ implies $0_W \in h(x_0)$,

Consequently, $\bar{x} \in D$ is a weakly efficient solution of (VP). □

Similar to the proof of Theorem 3.1, we have the following Theorem 3.2 and 3.3.

**Theorem 3.2** Let $\bar{x} \in D$. Assume that $f(x) - f(\bar{x}), g(x), h(x)$ satisfy the generalized convexity condition (a1) and the inner point condition (a2), and (VP) satisfies the Slater Constrained Qualification (SC). Then $\bar{x} \in D$ is a weakly efficient solution of (VP) if and only if $\exists \bar{y} \in f(\bar{x})$ and $\exists$Lagrangian multiplier $\xi \in Y_+^*, \eta \in Z_+^*, \varsigma \in W^*$ with $\xi \neq 0_{Y^*}$ such that

$$\min_{x \in D}[\xi(f(x)) + \eta(g(x)) + \varsigma(h(x))] \geq \xi(\bar{y}),$$
$$\min \eta(g(\bar{x})) = 0,$$

where, for the second equality, the minimum is taken for all $z \in g(\bar{x})$. □

**Theorem 3.3** Let $\bar{x} \in D$. And suppose $f(x) - f(\bar{x}), g(x), h(x)$ satisfy the generalized convexity condition (a1) and the inner point condition (a2), and (VP) satisfies the Slater Constraint Qualification (SC), then $\bar{x}$ is a weakly efficient solution of (VP) if and only if $\exists \xi \in Y_+^* \setminus \{0_{Y^*}\}$ such that $\bar{x}$ is an optimal solution of the following scalar optimization problem:

$$\min_{x \in D} \xi(f(x)). \square$$

Now, let $B(Z,Y)$ denote the set of continuous linear transformations from $Z$ to $Y$, and $B^+(Z,Y) = \{S \in B(Z,Y) : S(Z_+) \subseteq Y_+\}$.

**Definition 3.4** The vector Lagrangian map $L: X \times B^+(Z,Y) \times B(W,Y) \to 2^Y$ of (VP) is defined by the set-valued map

$$L(x, S, T) = f(x) + S(g(x)) + T(h(x)).$$

Given $(S, T) \in B^+(Z, Y) \times B(W, Y)$, we consider the minimization problem induced by (VP):

(VPST) $\quad\quad Y_+ - \min L(x, S, T),$
$\quad\quad\quad\quad\quad s.t., x \in D.$

According the following Theorem 3.4, (VPST) can also be considered as a dual problem of (VP).

**Theorem 3.4** Let $\bar{x} \in D$. Assume that $f(x) - f(\bar{x}), g(x), h(x)$ satisfy the generalized convexity condition (a1) and the inner point condition (a2), and (VP) satisfies the Slater Constrained Qualification (SC). Then $\bar{x} \in D$ is a weakly efficient solution of (VP) if and only if $\exists (S, T) \in B^+(Z, Y) \times B(W, Y)$ such that $\bar{x} \in D$ is a weakly efficient solution of (VPST).

**Proof**. Assume $\exists (S, T) \in B^+(Z, Y) \times B(W, Y)$ such that $\bar{x} \in D$ is a weakly efficient solution of (VPST). Then there exist $\bar{y} \in f(\bar{x}), \bar{z} \in g(\bar{x}), \bar{w} \in h(\bar{x})$, such that

$$(\bar{y} + S(\bar{z}) + T(\bar{w}) - [f(D) + S(g(D)) + T(h(D))]) \cap \text{int } Y_+ = \emptyset,$$

If $(\bar{y} - f(D)) \cap \text{int } Y_+ \neq \emptyset$, then $\exists y \in f(D)$ such that $\bar{y} - y \in \text{int } Y_+$, i.e.,

$$(\bar{y} + S(\bar{z}) + T(\bar{w}) - [y + S(\bar{z}) + T(\bar{w})]) \in \text{int } Y_+.$$

Which means that

$$(\bar{y} + S(\bar{z}) + T(\bar{w}) - [f(D) + S(g(D)) + T(h(D))]) \cap \text{int } Y_+ \neq \emptyset.$$

Which is a contradiction.

Therefore

$$(\bar{y} - f(D)) \cap \text{int } Y_+ = \emptyset.$$

Hence, from Lemma 3.1, $\bar{x} \in D$ is a weakly efficient solution of (VP).

Conversely, suppose that $\bar{x} \in D$ is a weakly efficient solution of (VP). So $\exists \bar{y} \in f(\bar{x})$ such that there is not any $x \in D$ for which $f(x) - \bar{y} \in -\text{int } Y_+$. That is to say, there is not any $x \in X$ such that

$$f(x) - \bar{y} \in -\text{int } Y_+, g(x) \in -Z_+, 0_W \in h(x).$$

By Theorem 2.2 $\exists (\xi, \eta, \varsigma) \in Y_+^* \times Z_+^* \times W^* \setminus \{(0_{Y^*}, 0_{Z^*}, 0_{W^*})\}$ such that

$$\xi(f(x) - \bar{y}) + \eta(g(x)) + \varsigma(h(x)) \geq 0, \forall x \in D. \tag{3.8}$$

Since $\bar{y} \in f(\bar{x})$ and $0_W \in h(\bar{x})$, take $x = \bar{x}$ in (3.8) we obtain

$$\eta(g(\bar{x})) \geq 0.$$

But $\bar{x} \in D$ and $\eta \in Z_+^*$ imply that $\exists \bar{z} \in g(\bar{x}) \cap (-Z_+)$ for which

$$\eta(\bar{z}) \leq 0.$$

Hence $\eta(\bar{z}) = 0$, which means

$$0 \in \eta(g(\bar{x})). \tag{3.9}$$

Since $x \in D$ implies $0_W \in h(x)$, and $g(x) \cap (-Z_+) \neq \emptyset$ implies $\exists z \in g(x) \cap (-Z_+)$ such that $\eta(z) \leq 0$, we have

$$\xi(f(x) - \bar{y}) \geq 0, \forall x \in D.$$

Because the Slater Constraint Qualification is satisfied, similar to the proof of Theorem 3.1, we have $\xi \neq 0_{Y^*}$. So we may take $y_0 \in \text{int } Y_+$ such that

$$\xi(y_0) = 1.$$

Define the operator $S : Z \to Y$ and $T : W \to Y$ by

$$S(z) = \eta(z)y_0, T(w) = \varsigma(w)y_0. \tag{3.10}$$

It is easy to see that

$$S \in B^+(Z,Y), S(Z_+) = \eta(Z_+)y_0 \subseteq Y_+,$$
$$T \in B(W,Y).$$

And (3.9) implies

$$S(g(\bar{x})) = \eta(g(\bar{x}))y_0 \in 0 \cdot Y_+ = 0_Y. \tag{3.11}$$

Since $\bar{x} \in D$, we have $0_W \in h(\bar{x})$. Hence

$$0_Y \in T(h(\bar{x})). \tag{3.12}$$

Therefore, by (3.11) and (3.12)

$$\bar{y} \in f(\bar{x}) \subseteq f(\bar{x}) + S(g(\bar{x})) + T(h(\bar{x})).$$

From (3.8) and (3.10)

$$\xi[f(x) + S(g(x)) + T(h(x))]$$
$$= \xi(f(x)) + \eta((g(x))\xi(y_0) + \varsigma(h(x))\xi(y_0)$$
$$= \xi(f(x)) + \eta(g(x)) + \varsigma(h(x))$$
$$\geq \xi(\bar{y}), \forall x \in D.$$

i.e.,

$$\xi[f(x) - \bar{y}) + S(g(x)) + T(h(x))]$$
$$\geq 0, \forall x \in D. \tag{3.13}$$

Taking $F(x) = f(x) + S(g(x)) + T(h(x))$, $G(x) = \{0_Z\}$ and $H(x) = \{0_W\}$, applying Theorem 2.2 to the functions $F(x) - \bar{y}, G(x), H(x)$, then (3.13) deduces that

$$(\bar{y} - [f(D) + S(g(D)) + T(h(D))]) \cap \text{int } Y_+ = \emptyset, \qquad (3.14)$$

and

$$\bar{y} \in F(\bar{x}) = f(\bar{x}) + S(g(\bar{x})) + T(h(\bar{x})),$$

since $0_Y \in S(g(\bar{x})), 0_Y \in T(h(\bar{x}))$.

Consequently, $\bar{x} \in D$ is a weakly efficient solution of (VPST).

We complete the proof. $\square$

Similar to the proof of Theorem 3.4, we have

**Theorem 3.5** Let $\bar{x} \in D$. Assume that $f(x) - f(\bar{x}), g(x), h(x)$ satisfy the generalized convexity condition (a1) and the inner point condition (a2). Then $\bar{x}$ is a weakly efficient solution of (VP) implies $\exists$vector Lagrangian multiplier $(S,T) \in B^+(Z,Y) \times B(W,Y)$ such that $\bar{x} \in D$ is a weakly efficient solution of (VPST).

Inversely, if (NNAMCQ) holds at $\bar{x} \in D$, and if $\exists$vector Lagrangian multiplier $(S,T) \in B^+(Z,Y) \times B(W,Y)$ such that $\bar{x}$ is a weakly efficient solution of (VPST), then $\bar{x}$ is a weakly efficient solution of (VP). $\square$

**Remark 3.1** We may have the generalized convexity conditions as follows:

$\exists$bounded set-valued maps $u_1 : X \to 2^Y, u_2 : X \to 2^Z$, $\forall \alpha \in (0,1), \forall x^i \in D$, $\forall y^i \in f(x^i), \forall z^i \in g(x^i), \forall w^i \in h(x^i), i = 1,2, \forall \varepsilon > 0$, such that

(b1)
$$\varepsilon u_1 + \alpha y^1 + (1-\alpha)y^2 \in f(D) + Y_+,$$
$$\varepsilon u_2 + \alpha z^1 + (1-\alpha)z^2 \in g(D) + Z_+,$$
$$\alpha w^1 + (1-\alpha)w^2 \in h(D).$$

Then, the condition "$f(x) - f(\bar{x}), g(x), h(x)$ satisfy the generalized convexity condition (a1)" in this paper can be replaced by "$f(x), g(x), h(x)$ satisfy the generalized convexity condition (b1)".